\documentclass{amsart}

\usepackage{amsfonts}
 \usepackage{graphicx}
  \usepackage{bbm}
   \usepackage{amscd}
    \usepackage{amsmath}

\usepackage{amssymb}
\usepackage{mathrsfs}
 \usepackage{epsfig} 

\newtheorem*{COP}{Theorem A}
\newtheorem*{mainthm}{Theorem B}
\newtheorem{theorem}{Theorem}[section]
\newtheorem{corollary}[theorem]{Corollary}
\newtheorem{lem}[theorem]{Lemma}

\theoremstyle{definition}

\theoremstyle{remark}
\newtheorem{rem}[theorem]{Remark}

\numberwithin{equation}{section}

\begin{document}

\title[A Positive Solution to a Conjecture of A.~Katok]
{A Positive Solution to a Conjecture of A.~Katok for Diffeomorphism Case}%
\author[X.~Dai]{Xiongping Dai}%
\address{Department of Mathematics\\ Nanjing University\\
Nanjing, 210093, P. R. CHINA}%
\email{xpdai@nju.edu.cn}%

\thanks{This project was supported by
NSFC (No.~10671088) and 973 project (No.~2006CB805903)}%

\subjclass[2000]{37D05, 37D20, 37D25}%
\keywords{Anosov diffeomorphism, Pesin stable manifold, Liao theory, shadowing property.}%

\date{June 9, 2008}
\begin{abstract}
A.~Katok has conjectured that a $C^{1+\alpha}$ map $g\colon
M^n\rightarrow M^n, n\ge 2$, which is H\"{o}lder conjugated to an
Anosov diffeomorphism is also an Anosov diffeomorphism. Using Pesin
stable manifold theorem and Liao spectrum theorem, we show that
under the hypothesis of such a conjecture, $g$ is an Axiom A
diffeomorphism having no cycles. Particularly, if $g$ is H\"{o}lder
conjugated to a hyperbolic toral automorphism, then $g$ is Anosov.
\end{abstract}
\maketitle

\section{\bf Introduction}\label{sec1}

Let $M^n$ be a connected, compact, smooth, and closed Riemannian
manifold of dimension $n\ge 2$. A.~Katok has conjectured that if
$g\in\mathrm{Diff}^{1+}(M)$ is H\"{o}lder conjugated to an Anosov
$f\in\mathrm{Diff}^{1}(M)$; i.e., there is a
H\"{o}lder-homeomorphism $h$ of $M$ such that $f=h\circ g\circ
h^{-1}$, then $g$ is also an Anosov diffeomorphism. Here
``H\"{o}lder-homeomorphism $h$" means that $h$ and its inverse
$h^{-1}$ both are H\"{o}lder continuous. And $\mathrm{Diff}^{1+}(M)$
is the set of all $C^1$ diffeomorphisms with $\alpha$-H\"{o}lder
derivatives for some H\"{o}lder exponent $\alpha$ with
$0<\alpha\le1$.

For convenience, if $g\in\textrm{Diff}^{1+}(M)$ is H\"{o}lder
conjugated to an Anosov diffeomorphism, then $g$ is temporally said
to be \emph{Katok}. Under the hypothesis of such a conjecture, in
\cite{COP} the authors proved the following.

\begin{COP}[\cite{COP}]
If $g$ is Katok, then all periodic points of $g$ have only non-zero
Lyapunov exponents, and such exponents are uniformly bounded away
from zero.
\end{COP}

In this paper, using Pesin stable manifold theorem, Liao spectrum
theorem and Liao reordering theorem, and shadowing property, based
on Theorem A above we obtain a positive solution to Katok's
conjecture as follows.

\begin{mainthm}
If $g\in\mathrm{Diff}^{1+}(M)$ is Katok, then $g$ is an Axiom A
diffeomorphism having no cycles.
\end{mainthm}

Consequently, if $g$ is H\"{o}lder conjugated to an Anosov
diffeomorphism $f$ that satisfies $\Omega(f)=M$ such as a hyperbolic
toral automorphism, then $g$ is also Anosov. In addition, if $g$ is
volume-preserving, then $g$ is Anosov too. Theorem~B shows that
Anosov diffeomorphisms have strong rigidity.

This paper is organized as follows. In $\S\ref{sec2}$ we will
introduce the Liao spectrum theorem and reordering theorem for $C^1$
differential systems on Euclidean spaces. Then we will prove a
semi-uniform ergodic theorem which provides us with a criterion from
nonuniform hyperbolicity to uniform hyperbolicity. In $\S\ref{sec3}$
we will prove an approximation theorem of ergodic measure by
periodic measures. In $\S\ref{sec4}$ we will first prove that a
Katok diffeomorphism is nonuniformly hyperbolic for any invariant
measures and then show that it is uniformly hyperbolic. We will
consider a volume-preserving Katok diffeomorphism in the last
section.

\section{\bf Liao spectrum and reordering theorem}\label{sec2}
In this section, we will introduce the Liao spectrum theorem and
Liao reordering theorem, which are basic in Liao theory. Then,
applying the Liao spectrum theorem, we will provide with a criterion
of uniform contraction.

For simplicity, let us consider throughout this section a
nonsingular autonomous system of $C^1$-differential equations in an
$(n+1)$-dimensional Euclidean $w$-space $\mathbb{E}^{n+1}, n\ge 2$
\begin{equation*}
\dot{w}=S(w)\quad w\in \mathbb{E}^{n+1},\ S(w)\in\mathbb
R^{n+1}-\{\textbf{0}\},
\end{equation*}
where we write $T_w\mathbb{E}^{n+1}=\mathbb R^{n+1}$ for all $w$ to
distinguish the $w$-state-space $\mathbb{E}^{n+1}$ from its tangent
$x$-space $\mathbb{R}^{n+1}$, which then naturally gives rise to a
$C^1$-flow on the state-space $\mathbb{E}^{n+1}$
\begin{equation*}
\phi\colon\mathbb R\times
\mathbb{E}^{n+1}\rightarrow\mathbb{E}^{n+1};\ (t,w)\mapsto t_\cdot
w.
\end{equation*}
It further induces, on the tangent bundle
$T\mathbb{E}^{n+1}=\mathbb{E}^{n+1}\times\mathbb R^{n+1}$, a smooth
linear skew-product flow
\begin{align*}
&\varPhi\colon\mathbb R\times\mathbb{E}^{n+1}\times\mathbb
R^{n+1}\rightarrow\mathbb{E}^{n+1}\times\mathbb R^{n+1};\
(t,(w,x))\mapsto(t_\cdot
w,{\varPhi_{t,w}}{x})\\
\intertext{where $\varPhi_{t,w}\colon \mathbb{R}^{n+1}\rightarrow
\mathbb{R}^{n+1};\  x\mapsto\frac{\partial\phi(t,w)}{\partial w} x$,
corresponding to the extended system}&\dot{w}=S(w), \quad
\dot{x}=S^\prime(w)x
\end{align*}
on the extended $(w,x)$-phase-space $\mathbb{E}^{n+1}\times
\mathbb{R}^{n+1}$.

\subsection{}\label{sec2.1}

Let $\mathbb{T}=\bigcup_{w\in\mathbb{E}^{n+1}}\mathbb{T}_w$, where
$\mathbb{T}_w=\left\{x\in\mathbb R^{n+1}\,|\,\langle
S(w),x\rangle=0\right\}$, denote the subbundle of the tangent bundle
$\mathbb E^{n+1}\times\mathbb R^{n+1}$ transversal to $S$ over
$\mathbb{E}^{n+1}$, called the \emph{transversal tangent bundle to
$S$ over $\mathbb{E}^{n+1}$}. Then there is another naturally
induced smooth linear skew-product flow
\begin{equation*}
\varPsi\colon \mathbb R\times \mathbb{T}\rightarrow \mathbb{T};\
(t,(w,x))\mapsto(t_\cdot w,\varPsi_{t,w}{x})
\end{equation*}
where along the fiber direction, $\varPsi_{t,w}\colon
\mathbb{T}_w\rightarrow \mathbb{T}_{t_\cdot w}$ is defined as the
component of $\varPhi_{t,w}x$ transversal to $S(t_\cdot w)$ for any
$(w,x)\in \mathbb{T}$; that is, $\varPhi_{t,w}x=rS(t_\cdot
w)+\varPsi_{t,w}{x}$, $\varPsi_{t,w}x\in\mathbb{T}_{t_\cdot w}$, for
some $r\in\mathbb{R}$. Particularly, let
$$
\mathscr
F_{1}^{*\sharp}=\left\{(w,x)\in\mathbb{T}\colon\|x\|=1\right\}
$$
be the unit transversal tangent bundle to $S$ over
$\mathbb{E}^{n+1}$. Then, there is a natural skew-product flow
\begin{equation*}
\varPsi^\sharp\colon \mathbb
R\times\mathscr{F}_1^{*\sharp}\rightarrow\mathscr{F}_1^{*\sharp};\
(t,(w,x))\mapsto(t_\cdot w,\varPsi_{t,w}^\sharp{x})
\end{equation*}
where
$$
\varPsi_{t,w}^\sharp{x}=\varPsi_{t,w}{x}/\|\varPsi_{t,w}^\sharp{x}\|.
$$

\subsection{}\label{sec2.2}

By the \emph{bundle of transversal orthogonal $n$-frames} over
$\mathbb{E}^{n+1}$, it means $\mathscr F_n^*$, where the fiber at
$w$ is given by
\begin{equation*}
\mathscr F_{n,w}^*=\left\{\gamma\in
\stackrel{n\textrm{-time}}{\overbrace{\mathbb{T}_w\setminus\{\textbf{0}\}\times\cdots\times\mathbb{T}_w\setminus\{\textbf{0}\}}}\,|\,\langle\textrm{col}_i{\gamma},\textrm{col}_j{\gamma}\rangle=0,\
1\le i\neq j\le n\right\}.
\end{equation*}
Here and in the future, for $1\le i\le n$
$$
\textrm{col}_i\colon(v_1,\ldots,v_n)\mapsto v_i.
$$
Using the well-known Gram-Schmidt orthogonalizing process, based on
$(\mathbb{E}^{n+1},\phi)$, we can obtain from $\varPsi$ the
well-defined skew-product flow on $\mathscr F_n^*$ as follows
\begin{equation*}
\chi^*\colon\mathbb R\times\mathscr F_n^*\rightarrow\mathscr F_n^*;\
(t,(w,\gamma))\mapsto(t_\cdot w, \chi_{t,w}^{*}\gamma).
\end{equation*}
The {\em bundle of transversal orthonormal $n$-frames} over
$\mathbb{E}^{n+1}$ is written as $\mathscr F_n^{*\sharp}$, where the
fiber at $w$ is defined as
\begin{equation*}
\mathscr F_{n,w}^{*\sharp}=\left\{{\gamma}\in \mathscr
F_n^*\colon\|\textrm{col}_j{\gamma}\|=1,\ j = 1, \ldots, n\right\}.
\end{equation*}
Furthermore, there is a natural skew-product flow based on
$(\mathbb{E}^{n+1},\phi)$
\begin{equation*}
\chi^{*\sharp}\colon\mathbb R\times\mathscr F_n^{*\sharp}\rightarrow
\mathscr F_n^{*\sharp};\ (t,(w,\gamma))\mapsto(t_\cdot w,
\chi_{t,w}^{*\sharp}\gamma)
\end{equation*}
called the \emph{Liao transversal orthonormal $n$-frame flow of
$S$}. For convenience, let
\begin{equation*}
\pi\colon\mathscr F_n^{*\sharp}\rightarrow\mathbb{E}^{n+1};\
(w,{\gamma})\mapsto w
\end{equation*}
be the bundle projection. Then, the following commutativity holds:
\begin{equation*}
t_\cdot
w=\phi(t,\pi(w,{\gamma}))=\pi(\chi^{*\sharp}(t,(w,\gamma)))\quad
\forall\,(t,(w,\gamma))\in\mathbb R\times\mathscr F_n^{*\sharp}.
\end{equation*}

\subsection{}\label{sec2.3}

Now, the so-called \emph{Liao qualitative functions} of $S$ are the
following
\begin{align*}
&\omega_i^*\colon\mathscr F_n^{*\sharp}\rightarrow {\mathbb R}\\
\intertext{for $i=1,\ldots,n$, given by}
&\omega_i^*(w,{\gamma})=\left.\frac{d}{dt}\right|_{t=0}\|\textrm{col}_i\circ{\chi_{t,w}^{*}}{\gamma}\|
\end{align*}
where $\chi_{t,w}^{*}\cdot\colon\mathscr
F_{n,w}^{*\sharp}\rightarrow \mathscr F_{n,t_\cdot w}^*$ is as in
$\S\ref{sec2.2}$ for any $(t,w)\in\mathbb{R}\times\mathbb{E}^{n+1}$.
Particularly, let
\begin{equation}\label{eq2.11}
\omega^*\colon\mathscr F_1^{*\sharp}\rightarrow {\mathbb R};\
(w,x)\mapsto\left.\frac{d}{dt}\right|_{t=0}\|\varPsi_{t,w}x\|.
\end{equation}
Since $S$ is of class $C^1$, these functions $\omega^*, \omega_i^*,
1\le i\le n$, all are well defined and continuous; see~\cite{L93,
Dai}.

\subsection{}\label{sec2.4}

Let ${\mathcal M}_{\textit{inv}}$ and ${\mathcal M}_{\textit{erg}}$
denote the set of all invariant Borel probability measures and
ergodic Borel probability measures of a dynamical system,
respectively. The following theorem equivalently describes the
Lyapunov characteristic spectrum of $(\mathbb{E}^{n+1},S)$.

\begin{theorem}[Liao spectrum theorem~\cite{L93,Dai}]\label{thm2.1}
Let $S(w), S^\prime(w)$ be bounded on $\mathbb{E}^{n+1}$ and assume
that there is an ergodic $\phi$-invariant Borel probability measure
$\mu$ on $\mathbb{E}^{n+1}$. Then, there exists a $\phi$-invariant
Borel subset $L(\mu)$ of $\mathbb{E}^{n+1}$ such that:
\begin{enumerate}
\item $\mu(L(\mu))=1$;

\item every point $w$ in $L(\mu)$ is Oseledets regular for $\varPsi$ based on $(\mathbb{E}^{n+1},\phi)$;

\item given any $P\in{\mathcal M}_{\textit{erg}}({\mathscr
F}_n^{*\sharp},\chi^{*\sharp})$ with marginal $\mu$; i.e.,
$\mu=P\circ\pi^{-1}$,
\begin{equation*}
{\rm\textsf{Sp}}_{\textsl{Lia}}^*(S,\mu):=\left\{\vartheta_i^*(P)\,|\,
i=1,\ldots,n\right\}
\end{equation*}
is just the Lyapunov spectrum of $\varPsi$ based on $(\phi,\mu)$,
counting with multiplicity and ignoring the order, which is called
the ``spectrum of transversal Lyapunov exponents" of $(S,\mu)$ and
is independent of the choices of $P$, where
\begin{equation*}
\vartheta_i^*(P):=\int_{{\mathscr
F}_n^{*\sharp}}\omega_i^*(w,{\gamma}) \,dP(w,{\gamma})\quad
\textrm{for } i=1,\ldots,n.
\end{equation*}
\end{enumerate}
\end{theorem}

\noindent Notice that under the hypotheses of Theorem~\ref{thm2.1}
above, there $P\in{\mathcal M}_{\textit{erg}}({\mathscr
F}_n^{*\sharp},\chi^{*\sharp})$ with marginal $\mu$ is always
existent from~\cite{Dai}.

\begin{theorem}[Liao reordering theorem~\cite{L93,Dai}]\label{thm2.2}
Under the hypotheses of Theorem~\ref{thm2.1}, let
${\rm\textsf{Sp}}_{\textsl{Lia}}^*(S,\mu)=\left\{\lambda_i^*(\mu)\,|\,
i=1,\ldots,n\right\}$ be the spectrum of transversal Lyapunov
exponents of $(S,\mu)$. If $i\mapsto\varrho(i)$ is any given
permutation of $\{1,\ldots,n\}$, then there is some
$P_\varrho\in{\mathcal M}_{\textit{erg}}({\mathscr
F}_n^{*\sharp},\chi^{*\sharp})$ with marginal $\mu$ such that
$$
\vartheta_i^*(P_\varrho)=\lambda_{\varrho(i)}^*(\mu)\quad\textrm{for
} i=1,\ldots,n.
$$
\end{theorem}

The above spectrum theorem~\ref{thm2.1} and reordering
theorem~\ref{thm2.2} will play an important role for the proof of
Theorem~B stated in $\S\ref{sec1}$.

\subsection{}\label{sec2.5}

Based on Theorem~\ref{thm2.1}, for any $P\in{\mathcal
M}_{\textit{erg}}({\mathscr F}_1^{*\sharp},\varPsi^{\sharp})$ with
marginal $\mu$ we have
$$
\vartheta^*(P):=\int_{{\mathscr F}_1^{*\sharp}}\omega^*(w,x)
\,dP(w,x)\in{\rm\textsf{Sp}}_{\textsl{Lia}}^*(S,\mu).
$$
Now, the following semi-uniform ergodic theorem will play an
important role for the proof of our main result.

\begin{theorem}\label{thm2.3}
Let $\Lambda$ be an $\phi$-invariant compact subset of
$\mathbb{E}^{n+1}$ and $\Delta$ an $\phi$-invariant Borel subset of
$\Lambda$ with total measure $1$; that is to say, $\mu(\Delta)=1$
for all $\mu\in\mathcal M_{\textit{inv}}(\Lambda,\phi|\Lambda)$. Let
$$
\mathbb{D}\colon\Delta\ni w\mapsto D(w)\subset\mathbb{T}_w
$$
be an $\mathfrak i$-dimensional $\varPsi$-invariant measurable
distribution for some integer $1\le\mathfrak i<n$. If
$\varPsi|\mathbb{D}$ has only negative Lyapunov exponents at almost
every $w\in\Delta$ and $\mathbb{D}$ is such that
$\lim_{\ell\to\infty}D(w_\ell)=D(w)$ provided, of course, that this
limit exists for $w_\ell\to w$ in $\Delta$, then
$\varPsi|\mathbb{D}$ is uniformly contracting.
\end{theorem}

\begin{proof}
Without any loss of generality, assume $\Lambda=\overline{\Delta}$.
Let
$$
\mathscr{F}_1^{*\sharp}(\Delta)=\left\{(w,x)\in\mathscr{F}_1^{*\sharp}\,|\,w\in\Delta
\textrm{ and }x\in D(w)\right\},
$$
which is $\varPsi^{\sharp}$-invariant. Then, it is easily seen that
$Y:=\overline{\mathscr{F}_1^{*\sharp}(\Delta)}$ is an
$\varPsi^{\sharp}$-invariant compact subbundle of
$\mathscr{F}_1^{*\sharp}$ over $\Lambda$, such that
$$
P(Y-\mathscr{F}_1^{*\sharp}(\Delta))=0\quad \textrm{for all
}P\in\mathcal M_{\textit{inv}}(Y,\varPsi^\sharp|Y).
$$
Moreover, according to Theorem~\ref{thm2.1} we obtain
$$
\int_Y\omega^*\,d P<0\quad\forall\,P\in\mathcal
M_{\textit{erg}}(Y,\varPsi^\sharp|Y),
$$
since $\varPsi|\mathbb{D}$ is nonuniformly contracting for any
$\mu\in\mathcal M_{\textit{erg}}(\Lambda,\phi|\Lambda)$. From the
continuous-time version of a semi-uniform theorem of \cite{SS}
(\cite[Lemma~3.1]{D07}), it follows that there exist constants
$\sigma>0$ and $T_0>0$ such that
\begin{align*}
&\int_Y\omega^*\,dP\le -\sigma\quad\forall\,P\in\mathcal
M_{\textit{erg}}(Y,\varPsi^\sharp|Y)\\
\intertext{and uniformly, for all $T\ge T_0$}
&\frac{1}{T}\int_0^T\omega^*(\varPsi^\sharp(t+s,y))\,dt\le -\sigma/2
\end{align*}
for all $y\in Y$ and for any $s\in\mathbb{R}$. Next, by the identity
$$
\frac{1}{T}\log\|\varPsi_{t,w}x\|=\frac{1}{T}\int_0^T\omega^*(\varPsi^\sharp(t,(w,x)))\,dt
$$
for all $(w,x)\in\mathscr{F}_1^{*\sharp}$ and for any $T\not=0$, we
can easily obtain that $\varPsi|\mathbb{D}$ is uniformly
contracting.

This proves the theorem.
\end{proof}

Notice here that if the distribution $\mathbb{D}$ is continuous;
that is, $D(w_\ell)\to D(w)$ provided that $w_\ell\to w$ in
$\Lambda$, then from~\cite{Cao} one can directly obtain the uniform
contraction of $\varPsi|\mathbb{D}$.
\section{\bf Shadowing property and approximation of ergodic measures}\label{sec3}

If $g\in\mathrm{Diff}^1(M)$ is conjugated to an Anosov
diffeomorphism, then $\overline{\textrm{Per}(g)}=\Omega(g)$ and $g$
has the shadowing property (see~\cite[Proposition~8.5]{S}).
Naturally, we ask if every ergodic measure of $g$ can be
approximated by periodic measures or not. Under our context, this is
the case.

Let $S$ be a $C^1$-differential system on $\mathbb{E}^{n+1}$ as in
$\S\ref{sec2}$ and $\Lambda$ an $\phi$-invariant nonempty compact
subset of $\mathbb{E}^{n+1}$. We say that $(\phi,\Lambda)$ has the
\emph{shadowing by periodic points property} provided that to any
$\epsilon>0$, there corresponds to some $\alpha>0$, such that for
any orbit arc $\phi([0,\tau],w)\subset\Lambda, \tau\ge2$ with
$\|w-\tau_\cdot w\|<\alpha$, there exists a periodic point
$p\in\Lambda$ with period $\tau$ satisfying $\|t_\cdot w-t_\cdot
p\|<\epsilon$ for all $t\in[0,\tau]$. In this case, we say that the
orbit $\phi(t,p)$ $\epsilon$-shadows the orbit arc
$\phi([0,\tau],w)$. Notice here that we require $p\in\Lambda$.

An $\phi$-invariant Borel probability measure $\mu$ on
$\mathbb{E}^{n+1}$ is called a \emph{periodic measure} if it is
supported on a periodic orbit of $\phi$; that is,
$\textrm{supp}(\mu)=\overline{\phi(\mathbb{R},p)}$ for some periodic
point $p$. The following result shows that periodic measures are
dense in $\mathcal {M}_{\textit{erg}}(\Lambda,\phi)$ under the weak
$*$-topology.

\begin{theorem}\label{thm3.1}
If the compact subsystem $(\Lambda, \phi)$ has the shadowing by
periodic points property, then periodic measures are dense in
$\mathcal {M}_{\textit{erg}}(\Lambda,\phi)$; that is, for any $\mu$
in $\mathcal {M}_{\textit{erg}}(\Lambda,\phi)$ there is a sequence
of periodic measures $(\mu_k)$ on $\Lambda$ such that $\mu_k\to\mu$
as $k$ tends to $\infty$.
\end{theorem}

\begin{proof}
Let $\mu\in\mathcal {M}_{\textit{erg}}(\Lambda,\phi)$ be
non-periodic, and let $Q_\mu(\Lambda,\phi)$ be the quasi-regular
point set of $(\Lambda,\mu,\phi)$; that is, $w\in
Q_\mu(\Lambda,\phi)$ if and only if
$$
\lim_{T\to\infty}T^{-1}\int_0^T\varphi(t_\cdot
w)\,dt=\int_\Lambda\varphi\,d\mu
$$
for all $\varphi\in C(\Lambda)$.

For any $w\in\Lambda$ and $T>0$, using the Riesz representation
theorem, we define the empirical measure $\mu_{w,T}$ on $\Lambda$ by
$$
\mu_{w,T}(\varphi)=\frac{1}{T}\int_0^T\varphi(t_\cdot
w)\,dt\quad\forall\,\varphi\in C(\Lambda).
$$
Then, given any Poisson stable (recurrent) point $\hat{w}\in
Q_\mu(\Lambda,\phi)$ we have $\mu_{\hat{w},T}\to\mu$ as $T\to\infty$
in the sense of weak $*$-topology; that is,
$\mu_{\hat{w},T}(\varphi)\to\mu(\varphi)$ for all $\varphi\in
C(\Lambda)$.

For any $\varphi\in C(\Lambda)$, let
$$\|\varphi\|_\infty={\sup}_{x\in\Lambda}\|\varphi(x)\|$$ and
$$
\|\varphi\|_\textrm{L}=\sup_{x,y\in\Lambda,
x\not=y}\frac{\|\varphi(x)-\varphi(y)\|}{\|x-y\|}.
$$
Then $\textrm{BL}(\Lambda)=\{\varphi\in
C(\Lambda);\,\|\varphi\|_\infty+\|\varphi\|_\textrm{L}<\infty\}$ is
dense in $(C(\Lambda),\|\cdot\|_\infty)$;
see~\cite[Theorem~11.2.4]{Dud}.

Now, by the shadowing by periodic points property and the recurrence
of the motion $\phi(t,\hat{w})$, we can choose a sequence of
periodic points $p_i$ with period $T_i\to\infty$ such that for any
$i$, $\|t_\cdot \hat{w}-t_\cdot p_i\|<1/i$ for all $t\in[0,T_i]$.
Then, it is easily seen that $\mu_{p_i}$, defined by
$$
\mu_{p_i}(\varphi)=\frac{1}{T_i}\int_0^{T_i}\varphi(t_\cdot
p_i)\,dt\quad\forall\,\varphi\in C(\Lambda),
$$
is an ergodic periodic measure of the subsystem $(\Lambda,\phi)$.

Next, for any $\varphi\in \textrm{BL}(\Lambda)$ we have
\begin{align*}
\lim_{i\to\infty}|\int\varphi\,d\mu-\int\varphi\,d\mu_{p_i}|&\le\lim_{i\to\infty}|\int\varphi\,d\mu-\int\varphi\,d\mu_{\hat{w},T_i}|\\
&{\quad}\quad+\limsup_{i\to\infty}|\int\varphi\,d\mu_{\hat{w},T_i}-\int\varphi\,d\mu_{p_i}|\\
&\le\limsup_{i\to\infty}\frac{1}{T_i}\int_0^{T_i}|\varphi(t_\cdot\hat{w})-\varphi(t_\cdot
p_i)|\,dt\\
&\le\limsup_{i\to\infty}\frac{1}{T_i}\int_0^{T_i}\|\varphi\|_\textrm{L}\|t_\cdot\hat{w}-t_\cdot
p_i\|\,dt\\
&\le\limsup_{i\to\infty}\frac{\|\varphi\|_\textrm{L}}{i}\\
&=0,
\end{align*}
which implies that $\mu_{p_i}\to \mu$ by the density of
$\textrm{BL}(\Lambda)$ in $(C(\Lambda),\|\cdot\|_\infty)$, as
required.

This proves the theorem.
\end{proof}

\section{\bf Hyperbolicity of Katok maps}\label{sec4}

In this section, we will finish the proof of our main result
Theorem~B stated in the Introduction, using the theorems introduced
in $\S\S~\ref{sec2}$ and \ref{sec3}.

\subsection{}Let $S$ be a $C^1$-differential system on
$\mathbb{E}^{n+1}$ as in $\S\ref{sec2}$ and $\Lambda$ an
$\phi$-invariant nonempty compact subset of $\mathbb{E}^{n+1}$.

\begin{theorem}\label{thm4.1}
Assume that the subsystem $(\phi,\Lambda)$ has the shadowing by
periodic points property, and, for each periodic point $p$ in
$\Lambda$, let $\lambda_1^*(p)\le\cdots\le\lambda_n^*(p)$ be the
spectrum of transversal Lyapunov exponents of $S$ at $p$, counting
with multiplicity. If there are some $\sigma<\varsigma$ such that
$\lambda_{1}^*(p)\le\sigma$ and $\lambda_{n}^*(p)\ge\varsigma$ for
all $p\in\mathrm{Per}(\Lambda,\phi)$, then for all $\mu\in\mathcal
{M}_{\textit{erg}}(\Lambda,\phi)$, $(S,\mu)$ has at least two
transversal Lyapunov exponents $\lambda_-^*(\mu)\le\sigma$ and
$\lambda_+^*(\mu)\ge\varsigma$.
\end{theorem}

\begin{proof}
Let $\mu\in\mathcal {M}_{\textit{erg}}(\Lambda,\phi)$ be
non-periodic. To prove the statement, it is enough to show that
$(S,\mu)$ has at least two transversal Lyapunov exponents
$\lambda_-^*(\mu)$ and $\lambda_+^*(\mu)$ such that
$\lambda_-^*(\mu)\le\sigma$ and $\lambda_+^*(\mu)\ge\varsigma$.

Let
$\mathscr{F}_n^{*\sharp}(\Lambda)=\left\{(w,\gamma)\in\mathscr{F}_n^{*\sharp}\,|\,w\in\Lambda\right\}$.
From Theorem~\ref{thm3.1}, we can take a sequence of periodic
measures, say $\{\mu_{p_i}\}$, in $\mathcal
{M}_{\textit{erg}}(\Lambda,\phi)$ with $\mu_{p_i}\to\mu$. By using
Theorems~\ref{thm2.1} and \ref{thm2.2}, we can choose some
$$
P_i\in\mathcal
{M}_{\textit{erg}}(\mathscr{F}_n^{*\sharp}(\Lambda),\chi^{*\sharp})\quad
\textrm{with marginal } \mu_{p_i}
$$
for all $i$, such that
$$
\lambda_1^*(p_i)=\vartheta_1^*(P_i)\le\cdots\le\vartheta_n^*(P_i)=\lambda_n^*(p_i).
$$
Since $\mathcal
{M}_{\textit{inv}}(\mathscr{F}_n^{*\sharp}(\Lambda),\chi^{*\sharp})$
is compact under the weak $*$-topology, there is no loss of
generality in assuming that $P_i\to P$ for some $P\in\mathcal
{M}_{\textit{inv}}(\mathscr{F}_n^{*\sharp}(\Lambda),\chi^{*\sharp})$
with marginal $\mu$; i.e., $P\circ\pi^{-1}=\mu$. Thus,
\begin{align*}
\lim_{i\to\infty}\int_{\mathscr{F}_n^{*\sharp}(\Lambda)}\omega_1^*\,dP_i=\int_{\mathscr{F}_n^{*\sharp}(\Lambda)}\omega_1^*\,dP\le\sigma\\
\intertext{and}
\lim_{i\to\infty}\int_{\mathscr{F}_n^{*\sharp}(\Lambda)}\omega_n^*\,dP_i=\int_{\mathscr{F}_n^{*\sharp}(\Lambda)}\omega_n^*\,dP\ge\varsigma
\end{align*}
for $\omega_1^*$ and $\omega_n^*$ both are continuous on
$\mathscr{F}_n^{*\sharp}(\Lambda)$. Then, by the classical ergodic
decomposition theorem we can choose at least two $P_-$ and $P_+$ in
$\mathcal
{M}_{\textit{erg}}(\mathscr{F}_n^{*\sharp}(\Lambda),\chi^{*\sharp})$
with marginal $\mu$ such that
$$
\lambda_-^*(\mu):=\int_{\mathscr{F}_n^{*\sharp}(\Lambda)}\omega_1^*\,dP_-\le\sigma
\quad\textrm{and}\quad\lambda_+^*(\mu):=\int_{\mathscr{F}_n^{*\sharp}(\Lambda)}\omega_n^*\,dP_+\ge\varsigma.
$$
By Theorem~\ref{thm2.1} again, $\lambda_-^*(\mu)$ and
$\lambda_+^*(\mu)$ both lie in
${\rm\textsf{Sp}}_{\textsl{Lia}}^*(S,\mu)$, as required.

This proves the theorem.
\end{proof}

\subsection{}To prove Theorem~B, we need a further remark on
Theorem A stated in $\S\ref{sec1}$. Let $g\in\mathrm{Diff}^{1+}(M)$
be H\"{o}lder conjugated to an Anosov $f\in\mathrm{Diff}^{1}(M)$.
Theorem~A asserts that all periodic points of $g$ have only non-zero
Lyapunov exponents. However, there is no information that if there
are any contracting or expanding periodic orbits. The following
lemma shows that a Katok map has no any contracting and expanding
periodic orbits.

Recall $f\in \textrm{Diff}^1(M)$ is said to be \emph{Anosov} if
there is a continuous splitting of $T_xM=E^s(x)\oplus E^u(x)$ for
every $x\in M$ and constants $C>0, \lambda>1$ such that
\begin{align*}
D_xf(E^s(x))&=E^s(f(x))\textrm{  and  }D_xf(E^u(x))=E^u(f(x)),\\
\|(D_xf^n)\vec{v}\|&\le
C\lambda^{-n}\|\vec{v}\|\quad\forall\,\vec{v}\in E^s(x),\
n\in\mathbb{N},\\
\|(D_xf^{-n})\vec{u}\|&\le
C\lambda^{-n}\|\vec{u}\|\quad\forall\,\vec{u}\in E^u(x),\
n\in\mathbb{N}.
\end{align*}
Then, the nonnegative integer
$$\textrm{Ind}(x):=\dim E^s(x)$$
for all $x\in M$ is called the \emph{index} of $f$ at $x$.

For any $f\in \textrm{Diff}^1(M)$ and $\delta>0$, as usual, for
$x\in M$ let
\begin{align*}
W^s(x)&=\left\{y\in M\,|\,\textrm{dist}(f^kx,f^ky)\to0\textrm{ as }
k\to\infty\right\}\\
\intertext{and} W_\delta^s(x)&=\left\{y\in
M\,|\,\textrm{dist}(f^kx,f^ky)\le\delta\textrm{ and }
\lim_{k\to\infty}\textrm{dist}(f^kx,f^ky)=0\right\}
\end{align*}
be the stable set and the local stable set of $f$ at $x$,
respectively. If $f$ is partially hyperbolic or $f$ is a
nonuniformly partially hyperbolic $C^{1+\alpha}$ diffeomorphism,
then $W^s(x)$ has local smooth manifold structure with
$T_xW^s(x)=E^s(x)$ a.e. (\cite{HPS,P76}).

Similarly, one can define the unstable manifolds $W^u(x)$ and
$W_\delta^u(x)$.

\begin{lem}\label{lem4.2}
If $g\in\mathrm{Diff}^{1+}(M)$ is H\"{o}lder conjugated to an Anosov
$f\in\mathrm{Diff}^{1}(M)$, then all periodic points of $g$ have
only non-zero Lyapunov exponents, and such exponents are uniformly
bounded away from zero, and
$$\mathrm{Ind}(p)\equiv\mathfrak i\quad\forall\,p\in\mathrm{Per}(g)$$
for some integer $\mathfrak i$ with $1\le\mathfrak i<n$.
\end{lem}

\begin{proof}
We first assert that for an Anosov $f\colon M^n\rightarrow M^n$,
there is an integer $\mathfrak i$ with $1\le\mathfrak i<n$ such that
$$
\textrm{Ind}(\hat{x},f)=\mathfrak i\quad \textrm{for all }\hat{x}\in
M.
$$
In fact, let $\Lambda_i=\{\hat{x}\in M\colon
\textrm{Ind}(\hat{x},f)=i\}$ for $i=0,1,\ldots,n$. Since the
splitting $T_{\hat{x}}M=E^s(\hat{x},f)\oplus E^u(\hat{x},f)$ is
continuous with respect to $\hat{x}\in M$, $\Lambda_i$ is closed and
further open in $M$ for all $i=0,1,\ldots,n$. Thus, every
$\Lambda_i$ is either equal to $\varnothing$ or to $M$. Clearly,
$\Lambda_0=\Lambda_n=\varnothing$. This shows the assertion.

Let $h\colon M\rightarrow M$ be a H\"{o}lder conjugacy from $g$ to
$f$. Given any $p\in\mathrm{Per}(g)$ and let $\hat{p}=h(p)$. From
Theorem~A it follows that $h(W^s(p;g))\subset W^s(\hat{p};f)$ and
$h(W^u(p;g))\subset W^u(\hat{p};f)$. Since $h$ is H\"{o}lder
homeomorphic, we have $\dim W_{\textit{loc}}^s(p;g)\le\dim
W_{\textit{loc}}^s(\hat{p};f)$ and $\dim
W_{\textit{loc}}^u(p;g)\le\dim W_{\textit{loc}}^u(\hat{p};f)$. Then,
we can obtain that $\dim W_{\textit{loc}}^s(p;g)=\dim
W_{\textit{loc}}^s(\hat{p};f)$ and so $\mathrm{Ind}(p)=\mathfrak i$
constant with $1\le\mathrm{Ind}(p)<n$.

This proves the lemma.
\end{proof}

\subsection{}Before proving Theorem~B, we first prove that a Katok map is non-uniformly hyperbolic.

\begin{theorem}\label{thm4.3}
Let $g\in\mathrm{Diff}^{1+}(M)$ be Katok. Then
\begin{enumerate}
\item there is a $\sigma>0$ such that to any $\mu\in\mathcal
{M}_{\textit{erg}}(M,g)$, $(g,\mu)$ is non-uniformly hyperbolic and
has at least two Lyapunov exponents, say $\lambda_-(\mu)$ and
$\lambda_+(\mu)$, with $\lambda_-(\mu)\le-\sigma$ and
$\lambda_+(\mu)\ge\sigma$;

\item there is an invariant subset $\varGamma$ and a measurable function $\delta\colon\varGamma\rightarrow(0,\infty)$ such
that $\mu\varGamma=1$ for all $\mu\in\mathcal
{M}_{\textit{erg}}(M,g)$ and $x\mapsto W_{\delta(x)}^s(x),\ x\mapsto
W_{\delta(x)}^u(x)$ both are well defined and continuous for $x$ in
$\varGamma$.
\end{enumerate}
\end{theorem}

\noindent Here $W_{\delta}^s(x)$ and $W_{\delta}^u(x)$ mean the
local stable and unstable manifolds of $g$ at $x$, respectively.

\begin{proof}
Let $g\in\mathrm{Diff}^{1+}(M)$ be H\"{o}lder conjugated to an
Anosov diffeomorphism $f\colon M\rightarrow M$. Given any
$\mu\in\mathcal {M}_{\textit{erg}}(M,g)$.

By the so-called suspension technique, from Lemma~\ref{lem4.2} and
Theorem~\ref{thm4.1} it follows immediately that $(g,\mu)$ has at
least two Lyapunov exponents, say $\lambda_-(\mu)$ and
$\lambda_+(\mu)$, such that $\lambda_-(\mu)\le-\sigma$ and
$\lambda_+(\mu)\ge\sigma$, where $\sigma$ is some positive constant
which is independent of $\mu$.

We next proceed to prove that $(g,\mu)$ is non-uniformly hyperbolic.
Let
$$T_xM=E^s(x,g)\oplus E^c(x,g)\oplus E^u(x,g)$$
for $\mu$-a.e.\,$x\in M$, where $E^s(x,g), E^c(x,g)$, and $E^u(x,g)$
stand for the stable direction, central direction and unstable
direction, respectively, associated to the Oseledets splitting of
$Dg$ at $x$~\cite{Os,L93}. We have $\dim E^s(x,g)\ge 1$ and $\dim
E^u(x,g)\ge 1$ for a.e. $x$. Since $g$ is of class $C^{1+\alpha}$
for some H\"{o}lder exponent $0<\alpha<1$, according to Pesin
theory~\cite{P76} there are local stable manifold
$W_{\textit{loc}}^s(x;g)$ and local unstable manifold
$W_{\textit{loc}}^u(x;g)$ with $\dim W_{\textit{loc}}^s(x;g)=\dim
E^s(x,g)$ and $\dim W_{\textit{loc}}^u(x;g)=\dim E^u(x,g)$ for
$\mu$-a.e. $x\in M$. On the other hand, by the
$C^\alpha$-conjugation $h^{-1}\colon M\rightarrow M$ from the Anosov
diffeomorphism $f$ to $g$, we obtain that for $\hat{x}=h(x)$
$$
\dim W_{\textit{loc}}^s(x;g)\ge\dim
W_{\textit{loc}}^s(\hat{x};f)\quad\textrm{and}\quad\dim
W_{\textit{loc}}^u(x;g)\ge\dim W_{\textit{loc}}^u(\hat{x};f),
$$
which implies $\dim W_{\textit{loc}}^s(x;g)+\dim
W_{\textit{loc}}^u(x;g)=n$ and so $E^c(x,g)=\textbf{0}$ for
$\mu$-a.e. $x\in M$. Thus, $(g,\mu)$ is non-uniformly hyperbolic.
This proves the statement (1).

Next, we are going to prove the statement (2). Since $f$ is Anosov,
there is some constant $\hat{\delta}>0$ such that the local stable
foliation $\mathscr{W}^s=(W_{\hat{\delta}}^s(\hat{x};f))_{\hat{x}\in
M}$ is continuous in $\hat{x}\in M$. Let $\varGamma$ be the
non-uniformly hyperbolic Pesin regular set of $g$~\cite{BP}.
Noticing that $h, h^{-1} $ both are H\"{o}lder and
$h^{-1}(W_{\hat{\delta}}^s(\hat{x};f))\subset W^s(x;g)$ where
$h(x)=\hat{x}$ for all $x\in\varGamma$, we can easily find some
measurable function $\delta\colon\varGamma\rightarrow(0,\infty)$,
which satisfies the requirements of the statement (2). This proves
the statement (2).

Thus, Theorem~\ref{thm4.3} is proved.
\end{proof}

\subsection{}
In~\cite{CLR}, the authors exhibit an example of a non-hyperbolic
horseshoe such that all Lyapunov exponents are non-zero and
uniformly bounded away from zero for all invariant measures. This
phenomenon is named ``completely nonuniformly hyperbolic."
Theorem~\ref{thm4.3} implies that a Katok map $g$ of $M^2$ has just
two Lyapunov exponents $\lambda_-(\mu)<0<\lambda_+(\mu)$, uniformly
bounded away from zero, for all $\mu$ in $\mathcal
{M}_{\textit{erg}}(M^2,g)$.

Nevertheless, there is still an essential gap from
Theorem~\ref{thm4.3} to Katok's conjecture even though in the
$2$-dimensional case, the continuity of the foliation
$(W_{\delta}^s(x))_{x\in\varGamma}$ guaranteed by
Theorem~\ref{thm4.3}(2), can avoid the occurrence of the completely
nonuniformly hyperbolic phenomenon.

Now we can finish the proof of Theorem~B using the semi-uniform
ergodic theorem Theorem~\ref{thm2.3}.

\begin{proof}[Proof of Theorem~B]
Let $\varGamma$ be defined by Theorem~\ref{thm4.3}(2) and let
$$T_xM=E^s(x)\oplus E^u(x)$$ for all $x\in\varGamma$ be the Oseledets
splitting of $g$ according to the multiplicative ergodic theorem.
Let $\mathfrak i$ be the index of $g$ and $\mathscr{G}_\mathfrak
i(T_{\varGamma}M)$ the Grassmannian manifold of $\mathfrak
i$-dimensional linear subspaces of $T_{\varGamma}M$. Then
$$
\mathbb{D}\colon\varGamma\rightarrow\mathscr{G}_\mathfrak
i(T_{\varGamma}M);\ x\mapsto E^s(x)
$$
is a $Dg$-invariant measurable distribution over $\varGamma$ such
that $T_xW_{\delta(x)}^s(x)=E^s(x)$ and $\dim E^s(x)=\mathfrak i$
for all $x\in\varGamma$.

Let $x_\ell\to x$ in $\varGamma$ and
$\lim_{\ell\to\infty}E^s(x_\ell)=E(x)$ for some $\mathfrak
i$-dimensional linear subspace $E(x)\subset T_xM$. As
$W_{\delta(x_\ell)}^s(x_\ell)\to W_{\delta(x)}^s(x)$ as
$\ell\to\infty$ by Theorem~\ref{thm4.3}, it follows from
$T_{x_\ell}W_{\delta(x_\ell)}^s(x_\ell)=E^s(x_\ell)$ that
$T_{x}W_{\delta(x)}^s(x)=E(x)$. Thus, $E(x)=E^s(x)$. Then, by the
discretization of Theorem~\ref{thm2.3} using the terms introduced
in~\cite{D05}, we obtain that $Dg\colon
\bigcup_{x\in\varGamma}E^s(x)\rightarrow\bigcup_{x\in\varGamma}E^s(x)$
is uniformly contracting.

Similarly, we can show that $Dg\colon
\bigcup_{x\in\varGamma}E^u(x)\rightarrow\bigcup_{x\in\varGamma}E^u(x)$
is uniformly expanding.

Therefore, $g$ is uniformly hyperbolic on $\overline{\varGamma}$.
Since $\Omega(g)=\overline{\textrm{Per}(g)}$, we have that
$\overline{\varGamma}=\Omega(g)$. Thus, $g$ is of Axiom A. Clearly,
$g$ has no cycles.

This proves the theorem.
\end{proof}

\section{\bf Volume-preserving Katok maps}\label{sec5}

Let $\textrm{Leb}$ denote the standard volume measure of $M^n$. A
Borel probability measure $\mu$ on $M$ is called a \emph{smooth
probability measure} if $\mu$ is absolutely continuous with respect
to $\textrm{Leb}$ such that
$$
C\le d\mu/d\textrm{Leb}\le K
$$
for some constants $C,K>0$. If $g$ is a $C^{1+\alpha}$
volume-preserving Anosov diffeomorphism, $g$ is ergodic due to
Anosov~\cite{Ano63,Ano}. For other proof of Anosov's ergodicity
theorem, see~\cite{X}. However, we are going to prove that if a
Katok map preserves a smooth probability measure, then it is Anosov
and thus also ergodic.

Using different approaches, it was proved independently by Bochi
$\&$ Viana~\cite{BV} and Xia~\cite{X} that the uniformly hyperbolic
closed sets of every $C^{1+\alpha}$ volume-preserving
diffeomorphisms have zero Lebesgue measure, unless they coincide
with the whole ambient compact manifold (Anosov case). This result
was generalized by using Pesin theory as follows:

\begin{lem}[\cite{D08}]\label{lem5.1}
Let $f$ be a diffeomorphism preserving a smooth probability measure
$\mu$ on a compact, connected, and closed Riemannian manifold $M$.
Let $\Lambda\subset M$ be a uniformly hyperbolic invariant Borel set
(not necessarily closed). If $\mu(\Lambda)>0$, then $f$ is Anosov
and $\Lambda=M\ (\textrm{mod }0)$.
\end{lem}

Now we can prove the ergodicity of a Katok map.

\begin{corollary}\label{cor5.2}
Let $g$ be a $C^{1+\alpha}$ Katok diffeomorphism of $M$ which
preserves a smooth probability measure $\mu$. Then $g$ is Anosov and
ergodic with $\Omega(g)=M$ .
\end{corollary}

\begin{proof}
Theorem~B and Lemma~\ref{lem5.1} imply that $g$ is Anosov with
$\Omega(g)=M$.

This proves Corollary~\ref{cor5.2}.
\end{proof}

\begin{rem}
We noted that it was recently announced by Zhihong Xia~\cite{SX}
that an Anosov diffeomorphism must be topologically transitive.
Then, our Theorem~B implies that every Katok diffeomorphism must be
Anosov.
\end{rem}
\subsection*{Acknowledgments}The author is grateful to Professors Huyi Hu and
Yunping Jiang for discussions on this paper.


\begin{thebibliography}{99}
\bibitem{Ano63}D.\,V.~Anosov, Ergodic properties of geodesic flows on closed Riemannian
manifolds of negative curvature. \textit{Soviet Math. Dokl.},
\textbf{4} (1963), 1153--1156.
\bibitem{Ano}D.\,V.~Anosov, Geodesic flows on closed Riemannian
manifolds of negative curvature. \textit{Proc. Steklov Math. Inst.},
\textbf{90} (1967), 1--235.

\bibitem{BP}L.~Barreira and Ya.~Pesin, Lectures on Lyapunov Exponents and Smooth Ergodic Theory. Appendix A by M. Brin and
Appendix B by D. Dolgopyat, H. Hu and Pesin. \textit{Proc. Sympos.
Pure Math.}, \textbf{69}, Smooth ergodic theory and its applications
(Seattle, WA, 1999), 3--106, Amer. Math. Soc., Providence, RI, 2001.

\bibitem{BV}J.~Bochi and M.~Viana, Lyapunov exponents: How
frequently are dynamical systems hyperbolic? in \textit{Modern
Dynamical Systems and Applications}, 271--297, Cambridge Univ.
Press, Cambridge, 2004.

\bibitem{Cao}Y.~Cao, Non-zero Lyapunov exponents and uniform
hyperbolicity. \textit{Nonlinearity}, \textbf{16} (2003),
1473--1479.

\bibitem{CLR}Y.~Cao, S.~Luzzato and I.~Rios, Some non-hyperbolic
systems with strictly non-zero Lyapunov exponents for all invariant
measures: horseshoes with internal tangencies. \textit{Discrete
Contin. Dyn. Syst.}, \textbf{15} (2006), 61--71.

\bibitem{COP}A.~Castro, K.~Oliveira and V.~Pinheiro, Shadowing by
non-uniformly hyperbolic periodic points and uniform hyperbolicity.
\textit{Nonlinearity}, \textbf{20} (2007), 75--85.

\bibitem{D05}X.~Dai, Partial linearization of differentiable systems. \textit{J. Difference Equ. Appl.}, \textbf{11} (2005), 965--977.
\bibitem{D07}X.~Dai, Hyperbolicity and integral expression of the
Lyapunov exponents for linear cocycles. \textit{J. Differential
Equations}, \textbf{242} (2007), 121--170.
\bibitem{Dai}X.~Dai, Integral expressions of Lyapunov exponents for
autonomous ordinary differential systems. \textit{Science in China
Series A: Mathematics}, \textbf{51} (2008), 000--000.
\bibitem{D08}X.~Dai, $C^{1+\alpha}$ volume preserving
diffeomorphisms have no any fat hyperbolic sets. Preprint 2008.

\bibitem{Dud}R.\,M.~Dudley, \textit{Real Analysis and Probability}. Second edition. Cambridge University Press, Cambridge,
2003.

\bibitem{HPS}M.~Hirsch, C.~Pugh and M.~Shub, Invariant manifolds.
\textit{Lect. Notes in Math.}, \textbf{583}, Springer-Verlag, 1977.

\bibitem{L93}S.-T.~Liao, On characteristic exponents construction of
a new Borel set for the multiplicative ergodic theorem for vector
fields. \textit{Acta Sci. Natur. Univ. Pekinensis}, \textbf{29}
(1993), 277--302.

\bibitem{Os}V.\,I.~Oseledec, A multiplicative ergodic theorem, Lyapunov characteristic numbers for dynamical systems.
\textit{Trudy Mosk Mat Obsec}, {\bf 19} (1968), 119--210.

\bibitem{P76}Ya.~Pesin, Families of invariant manifolds
corresponding to nonzero characteristic exponents, \textit{Math.
USSR-Izv.} \textbf{10} (1976), 1261--1305.

\bibitem{S}M.~Shub, \textit{Global Stability of Dynamical Systems}. Springer-Verlag, Berlin and
Heidelberg, 1987.

\bibitem{SX}R.~Saghin and Z. Xia, Homology of invariant foliations and its applications in
dynamics. International Conference on Topology and its Applications,
December 3--7, 2007 Kyoto, Japan.

\bibitem{SS}R.~Sturman and J.~Stark, Semi-uniform ergodic theorems
and applications to forced systems. \textit{Nonlinearity},
\textbf{13} (2000), 113--143.

\bibitem{X}Z.~Xia, Hyperbolic invariant sets with positive measures. \textit{Discrete
Contin. Dyn. Syst.}, \textbf{15} (2006), 811--818.
\end{thebibliography}
\end{document}